# Clapp-Puppe Type Lusternik-Schnirelmann (Co)Category in a Model Category

Donald Yau

Abstract. We introduce Clapp-Puppe type generalized Lusternik-Schnirelmann (co)category in a Quillen model category. We establish some of their basic properties and give various characterizations of them. As the first application of these characterizations, we show that our generalized (co)category is invariant under Quillen modelization equivalences. In particular, generalized (co)category of spaces and simplicial sets coincide. Another application of these characterizations is to define and study rational cocategory. Various other applications are also given.

## Introduction

Clapp and Puppe [**5, 6**] generalized the classical Lusternik-Schnirelmann category by introducing a class $\mathcal{A}$ of spaces and called the new invariant $\mathcal{A}$-category. The usual LS category of spaces is obtained by taking $\mathcal{A} = \{*\}$ (or the class of contractible spaces). In the 1960s, Ganea [**15, 16, 17**] introduced an invariant for spaces which he called cocategory as a dual to LS category within the framework of Eckmann-Hilton duality. Related notions of cocategory were later introduced by, in chronological order, Varadarajan [**32**], Hopkins [**23**], Hovey [**24**], and Doeraene [**9, 10**]. See James' survey article [**26**] for general information about LS category.

The purpose of this paper is to study generalized (co)category (in the sense of Clapp-Puppe) in a Quillen model category. On the one hand, our generalization of (co)category, denoted $\mathcal{A}\,\mathrm{cat}_{\mathbf{C}}$ ($\mathcal{A}\,\mathrm{cocat}_{\mathbf{C}}$), when applied to pointed spaces, is just Clapp-Puppe's $\mathcal{A}$-category, and if we further take $\mathcal{A}$ to be the class of contractible spaces, then we recover Ganea's (co)category. On the other hand, if $\mathcal{A}$ consists of the contractible objects in a model category $\mathbf{C}$, then our $\mathcal{A}\,\mathrm{cat}_{\mathbf{C}}$ reduces to Doeraene's $\mathbf{C}\,\mathrm{cat}$. So our work can be regarded as a common generalization of those of Clapp-Puppe and of Doeraene.

This paper is organized as follows. In §1 we introduce our Clapp-Puppe type generalization of LS (co)category in a Quillen model category and establish some basic properties. The main purpose of §2 is to obtain various characterizations of our generalized cocategory (resp. category) in terms of homotopy pullbacks (resp. pushouts). In §3.1 we will use these characterizations to prove the invariance of generalized (co)category under Quillen modelization equivalences. Then in §3.2 we define rational cocategory in terms of Sullivan's minimal model and use the results in §2 to prove that Ganea's cocategory rationalizes to our rational cocategory and that our rational cocategory satisfies the (rational) fibration property. Section 3.3 contains formulae for the generalized cocategory (resp. category) of a finite Cartesian product (resp. coproduct) and the generalized cocategory of a homotopy function complex.

We assume that the reader is familiar with the language of model category, of which [**12, 25**] are good references.

I would like to express my deepest gratitude to Professor Haynes Miller for many hours of inspiring conversations and encouragement. I would also like to thank the referee for helpful suggestions.

Contents

## 1. $\mathcal{A}$-category and $\mathcal{A}$-cocategory

Throughout this paper, we work with pointed proper model categories, usually denoted by $\mathbf{C}$, $\mathbf{D}$ etc. For an object $X$ in a model category, $QRX$ denotes, as usual, its functorial fibrant-cofibrant replacement. $\mathbf{C}_{cf}$ denotes the full subcategory of $\mathbf{C}$ consisting of all objects that are both fibrant and cofibrant. We call these objects fibrant-cofibrant. The symbols $\twoheadrightarrow$, $\rightarrowtail$, and $\xrightarrow{\sim}$ denote, respectively, a fibration, cofibration, and weak equivalence.

The purposes of this section are to introduce $\mathcal{A}\,\mathrm{cat}_{\mathbf{C}}$ and $\mathcal{A}\,\mathrm{cocat}_{\mathbf{C}}$ and to establish some basic properties, which will set the stage for our work in §2 and §3.

In §1.1 we give a brief account of Clapp-Puppe's generalization of LS category. Our definition of $\mathcal{A}$-cocategory is given in §1.2, in which we also show that our $\mathcal{A}\,\mathrm{cocat}$ is well-defined and is invariant under weak equivalences. It is shown in §1.3 that there is no loss of generality in assuming that $\mathcal{A}$ is closed under weak equivalence (within $\mathbf{C}_{cf}$). In §1.4 we show that if an object $Y$ is a retract of another object $X$, then the $\mathcal{A}$-cocategory of $Y$ is no bigger than that of $X$. In §1.5 we present the dual definition of $\mathcal{A}$-category and corresponding dual results for $\mathcal{A}\,\mathrm{cat}$.

**1.1. Clapp-Puppe's generalization of LS category.** Our definition of generalized (co)category is based on the generalized category introduced by Clapp and Puppe [**5, 6**]. In this subsection we give a summary of some of their main results.

Let $\mathcal{A}$ be a class of spaces with at least one non-empty space. The *$\mathcal{A}$-category* of a space $X$, denoted $\mathcal{A}\,\mathrm{cat}(X)$, is the smallest integer $n \geq 0$ such that there exists an open covering $\{X_0, \dots, X_n\}$ of $X$ and each inclusion $X_i \hookrightarrow X$ factors through some space in $\mathcal{A}$ up to homotopy; if no such integer $n$ exists we put $\mathcal{A}\,\mathrm{cat}(X) = \infty$. Such a covering is said to be *$\mathcal{A}$-categorical.* The *strong $\mathcal{A}$-category* $\mathcal{A}\,\mathrm{Cat}(X)$ of $X$ is the least integer $n \geq 0$ (or $\infty$ if no such $n$ exists) for which $X$ is homotopy equivalent to a space $X'$ that has an open covering $\{X'_0, \dots, X'_n\}$ such that each $X'_i$ has the homotopy type of some space in $\mathcal{A}$. $\mathcal{A}\,\mathrm{cat}(X)$ and $\mathcal{A}\,\mathrm{Cat}(X)$ are invariants of the homotopy type of $X$, and the former satisfies a general version of the Eilenberg cup-length theorem.

For many purposes it is convenient to have a certain "universal" map $f\colon U \to X$ from a space $U \in \mathcal{A}$ into $X$ with the property that every map from a space $A \in \mathcal{A}$ into $X$ factors through $f$ (not necessarily uniquely) up to homotopy. Such a map $f$ is called an *$\mathcal{A}$-left-versal* map for $X$. (Clapp and Puppe called it *$\mathcal{A}$-universal.*) The reader is referred to [**5**, 4.2] for examples of classes $\mathcal{A}$ for which every space has a left-versal map.

Let $f\colon Y \to X$ be a map. Define a sequence of fibrations $f_n\colon Y_n \to X$ inductively as follows. Let $f_1$ be the fibration associated with $f$. Suppose $f_k$ is defined for some $k \geq 1$. Form the pullback (which is also a homotopy pullback) $Y_1 \times_X Y_k = \lim \left(Y_1 \xrightarrow{f_1} X \xleftarrow{f_k} Y_k\right)$, let $Y'_{k+1}$ be the homotopy colimit of the diagram $(Y_1 \leftarrow Y_1 \times_X Y_k \rightarrow Y_k)$, and let $f_{k+1}\colon Y_{k+1} \to X$ be the fibration associated with the map that is induced by $f_1$ and $f_k$.

It is clear that the classical construction of the Milnor filtration $\{B_n\Omega X\}$ is (up to homotopy) the case when $Y = *$.

Now given a map $f\colon Y \to X$, replace it by its mapping cylinder $f'\colon Y \to M(f)$, identify $Y$ with the top of the mapping cylinder, and define the *$k$-fold* *fat-wedge* $W_k^Y X$ of $X$ under $Y$ to be the subspace of the Cartesian product $M(f)^k$ consisting of all points with at least one coordinate in $Y$.

THEOREM 1.1 ([**5**] 4.6, 4.8, 4.9, 5.5). *Let* $f\colon Y \to X$ *be an* $\mathcal{A}$*-left-versal map for* $X$*. Then the following are equivalent:*

1. $\mathcal{A}\,\mathrm{cat}(X) \leq k$*;*
2. *the fibration* $f_{k+1}$ *has a section;*
3. *the diagonal map* $\Delta\colon M(f) \to M(f)^{k+1}$ *factors through the* $(k+1)$*-fold fat-wedge* $W_{k+1}^Y X$ *of* $X$ *under* $Y$ *up to homotopy;*
4. $X$ *is a homotopy retract of a space* $Z$ *with* $\mathcal{A}\,\mathrm{Cat}(Z) \leq k$*.*

THEOREM 1.2 ([**6**] Thm. 1, 7, 7′). *Let* $\mathcal{A}$*,* $\mathcal{B}$ *and* $\mathcal{C}$ *be classes of spaces satisfying:*

1. $A \times B \in h\mathcal{C}$ *for all spaces* $A \in \mathcal{A}$ *and* $B \in \mathcal{B}$*;*
2. $h\mathcal{C}$ *is closed under finite disjoint unions.*

*Then for any spaces* $X$ *and* $Y$*,*

$$\mathcal{C}\,\mathrm{cat}(X \times Y) \leq \mathcal{A}\,\mathrm{cat}(X) + \mathcal{B}\,\mathrm{cat}(Y).$$

*The same inequality holds with* Cat *replacing* cat *throughout.*

**1.2.** $\mathcal{A}$**-cocategory.** In this subsection we give the definition of the $\mathcal{A}$-cocategory of an object $X$ and show that it is invariant under weak equivalences.

First we need some definitions.

DEFINITION 1.3. Let $\mathcal{A} \subseteq \mathbf{C}_{cf}$ be a non-empty class of fibrant-cofibrant objects in $\mathbf{C}$ and let $X \in \mathbf{C}_{cf}$. Let $f\colon X \to U \in \mathcal{A}$ be a map such that for every map $g\colon X \to A \in \mathcal{A}$, there exists a (not necessarily unique) map $h\colon U \to A$ such that $g \sim h \circ f$. Such a map $f$ is said to be an $\mathcal{A}$-*right-versal map for* $X$. The class $\mathcal{A}$ is said to be *right-versal* if every fibrant-cofibrant object in $\mathbf{C}$ has an $\mathcal{A}$-right-versal map.

Let $\mathcal{A}$ be a right-versal class in $\mathbf{C}$. We want to define the $\mathcal{A}$-cocategory, $\mathcal{A}\,\mathrm{cocat}_{\mathbf{C}}(X)$, of a fibrant-cofibrant object $X$ in $\mathbf{C}$. Let $f\colon X \to U$ be an $\mathcal{A}$-right-versal map for $X$. We construct a sequence of cofibrations $f_n\colon X \to U_n$ as follows. The map $f_1$ is the canonical cofibration associated to $f$, i.e., apply the functorial factorization to $f = pf_1$ so that $f_1$ is a cofibration and $p$ is an acyclic fibration. Suppose the cofibration $f_k\colon X \to U_k$ is defined for some $k \geq 1$. Take the (homotopy) pushout $U_1 \vee_X U_k = \mathrm{colim}\left(U_1 \xleftarrow{f_1} X \xrightarrow{f_k} U_k\right)$, apply functorial factorization (an acyclic cofibration followed by a fibration) to the induced map $\bar{f}_k \colon U_1 \to U_1 \vee_X U_k$ and take the pullback:

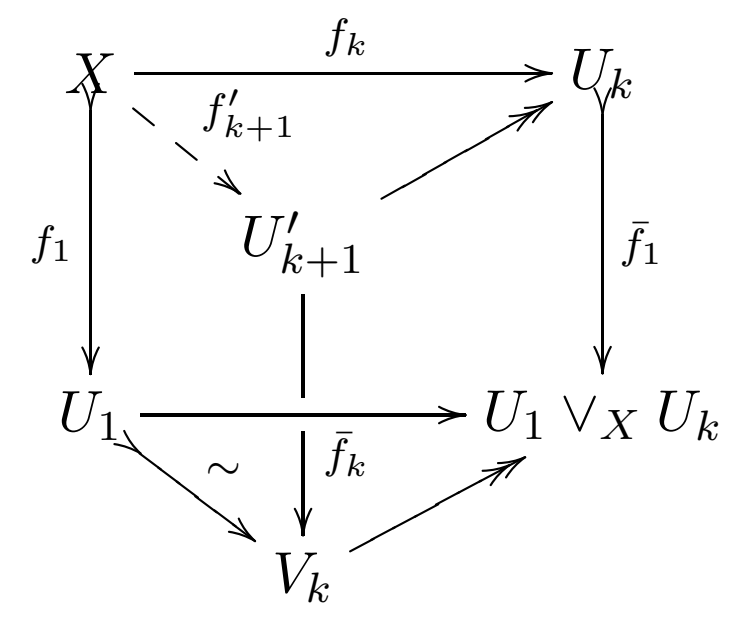

where we set $U''_k = U_1 \vee_X U_k$ and $U'_{k+1} = V_k \times_{U''_k} U_k$. Take $f_{k+1} \colon X \to U_{k+1}$ to be the canonical cofibration associated with the induced map $f'_{k+1} \colon X \to U'_{k+1}$.

Definition 1.4. Let $\mathcal{A}$ be a right-versal class in $\mathbf{C}$. The $\mathcal{A}$-*cocategory* of a fibrant-cofibrant object $X$, denoted $\mathcal{A}\,\mathrm{cocat}_{\mathbf{C}}(X)$ (or just $\mathcal{A}\,\mathrm{cocat}(X)$ when there is no danger of confusion), is the least integer $n \geq 0$ (possibly $\infty$) such that $f_{n+1}$ has a homotopy retraction, i.e., a map $r \colon U_{n+1} \to X$ such that $rf_{n+1} \sim \mathrm{id}_X$.

Whenever the symbol $\mathcal{A}\,\mathrm{cocat}_{\mathbf{C}}$ (or $\mathcal{A}\,\mathrm{cocat}$) occurs it is understood that $\mathcal{A}$ is a right-versal class. The definition of the $\mathcal{A}$-cocategory of an arbitrary object $X$ is given in Definition 1.9 below.

Remark 1.5. In the above construction each $U_k$ is both fibrant and cofibrant. Also, each $U'_k$ is fibrant because there is a fibration from it to $U_k$.

Proposition 1.6. *Let $X$ be a fibrant-cofibrant object in $\mathbf{C}$. Then $\mathcal{A}\,\mathrm{cocat}_{\mathbf{C}}(X)$ is independent of the choice of an $\mathcal{A}$-right-versal map for $X$.*

This proposition is a straightforward consequence of the following lemma.

Lemma 1.7. *If $f \colon X \to U$ and $g \colon X \to S$ are both $\mathcal{A}$-right-versal maps for $X$, then for each integer $k \geq 1$ there exists a map $\beta_k \colon U_k \to S_k$ in $\mathbf{C}$ such that $\beta_k f_k \sim g_k$.*

For the proof of Lemma 1.7 we will use the following standard terminology. Given a diagram $\{X \xrightarrow{f} Z \xleftarrow{g} Y\}$ the induced map from the pullback of the diagram, $X \times_Z Y$, to $Y$ is denoted by $f^\sharp$ and is called the *base change* of $f$ (along $g$). There is also a dual notion of *cobase change*.

Proof of Lemma 1.7. The proof is by induction. The $\mathcal{A}$-right-versality of $f$ implies that there exists a map $\beta'_1 \colon U \to S$ in $\mathbf{C}$ such that $\beta'_1 f \sim g$. So there is a homotopy commutative diagram in $\mathbf{C}$:

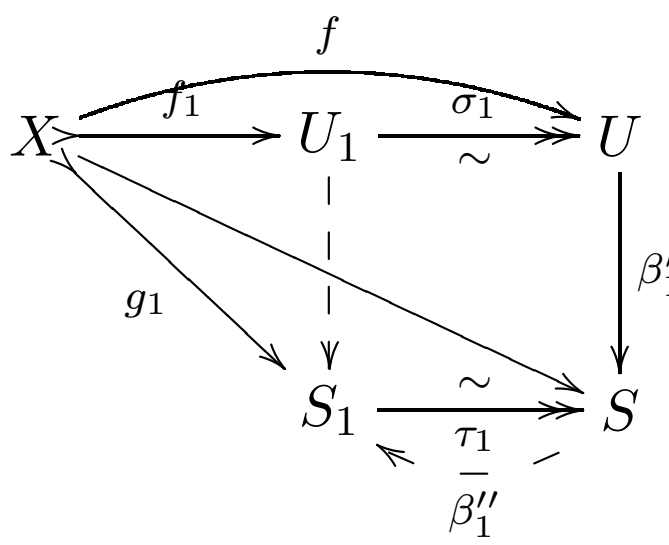

By [**25**, 1.2.8] a map $h\colon X_1 \to X_2$ between fibrant-cofibrant objects is a weak equivalence if and only if it is a homotopy equivalence. So there exists a homotopy inverse $\beta_1''\colon S \to S_1$ such that $\beta_1''\tau_1 \sim \mathrm{id}_{S_1}$. Now it follows that

$$g_1 \sim \beta_1''\tau_1 g_1 = \beta_1'' g \sim \beta_1''\beta_1' f = \beta_1''\beta_1'\sigma_1 f_1.$$

So we can take $\beta_1$ to be $\beta_1''\beta_1'\sigma_1$. This gets the induction started.

Suppose that for some $k \geq 1$ the map $\beta_k\colon U_k \to S_k$ in $\mathbf{C}$ exists such that $\beta_k f_k \sim g_k$. Consider the following commutative diagram in $\operatorname{Ho}\mathbf{C}$:

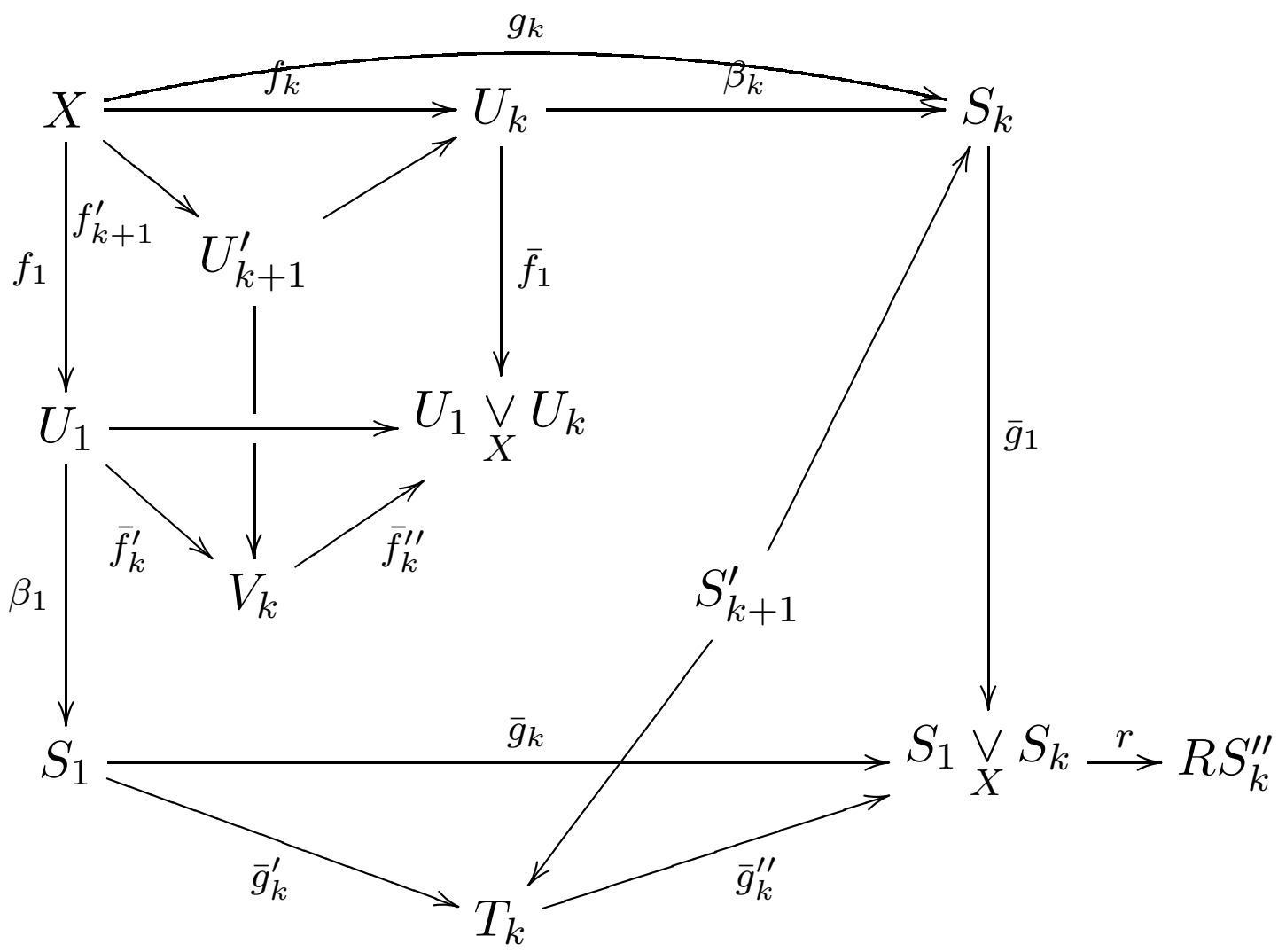

in which $U_{k+1}' = V_k \times_{U_k''} U_k$, $S_{k+1}' = T_k \times_{S_k''} S_k$, where $U_k'' = U_1 \vee_X U_k$ and $S_k'' = S_1 \vee_X S_k$ are the (homotopy) pushouts. The map $r$ is the functorial fibrant replacement of $S_k''$. By the results in [**12**, §10] the maps $r\bar{g}_k\beta_1\colon U_1 \to S_k''$ and $r\bar{g}_1\beta_k\colon U_k \to S_k''$ induce a map $\delta\colon U_k'' \to S_k''$ in $\operatorname{Ho}\mathbf{C}$ such that $\delta\bar{f}_k = \bar{g}_k\beta_1$ and $\delta\bar{f}_1 = \bar{g}_1\beta_k$. The map $\bar{f}_k'\colon U_1 \to V_k$ is an acyclic cofibration in $\mathbf{C}$ and so is an isomorphism in $\operatorname{Ho}\mathbf{C}$. An argument similar to the one above, using the universal property of $S_{k+1}'$, yields a map $\beta_{k+1}'\colon U_{k+1}' \to S_{k+1}'$ in $\operatorname{Ho}\mathbf{C}$ such that $\bar{g}_1^\sharp\beta_{k+1}' = \bar{g}_k'\beta_1(\bar{f}_k')^{-1}\bar{f}_1^\sharp$ and $\bar{g}_k''^\sharp\beta_{k+1}' = \beta_k\bar{f}_k''^\sharp$. A diagram-chase now shows that in $\operatorname{Ho}\mathbf{C}$ the composite $\beta_{k+1}'f_{k+1}'\colon X \to S_{k+1}'$ satisfies $\bar{g}_1^\sharp\beta_{k+1}'f_{k+1}' = \bar{g}_k'\beta_1(\bar{f}_k')^{-1}\bar{f}_1^\sharp f_{k+1}' = \bar{g}_k'\beta_1 f_1 =$ $\bar{g}'_k g_1$ and $\bar{g}''^{\sharp}_k \beta'_{k+1} f'_{k+1} = \beta_k \bar{f}''^{\sharp}_k f'_{k+1} = \beta_k f_k = g_k$. Thus it follows from the uniqueness of $g'_{k+1}\colon X \to S'_{k+1}$ that $g'_{k+1} = \beta'_{k+1} f'_{k+1}$ in $\operatorname{Ho}\mathbf{C}$. Now write $f'_{k+1} = f''_{k+1} f_{k+1}$ and $g'_{k+1} = g''_{k+1} g_{k+1}$, each as a cofibration followed by an acyclic fibration. The composite $\beta_{k+1} := (g''_{k+1})^{-1} \beta'_{k+1} f''_{k+1}$ in $\operatorname{Ho}\mathbf{C}$ satisfies $\beta_{k+1} f_{k+1} = (g''_{k+1})^{-1} \beta'_{k+1} f'_{k+1} = (g''_{k+1})^{-1} g'_{k+1} = (g''_{k+1})^{-1} g''_{k+1} g_{k+1} = g_{k+1}$. But since both $U_{k+1}$ and $S_{k+1}$ are fibrant and cofibrant (see Remark 1.5), the map $\beta_{k+1}$ in $\operatorname{Ho}\mathbf{C}$ lifts to a map (which we also call) $\beta_{k+1}$ in $\mathbf{C}$ that has the required property. This finishes the induction and also the proof of Lemma 1.7. □

Definition 1.8. Let $X$ be an object in a model category $\mathbf{C}$. By a *fibrant-cofibrant replacement* of $X$ we mean a fibrant-cofibrant object $W$ in $\mathbf{C}$, together with a zig-zag of weak equivalences connecting $X$ with $W$.

Definition 1.9. Let $\mathbf{C}$ be a pointed proper model category, let $\mathcal{A}$ be a right-versal class in $\mathbf{C}$, and let $X$ be an arbitrary object in $\mathbf{C}$. Then the $\mathcal{A}$-*cocategory* of $X$, denoted $\mathcal{A}\operatorname{cocat}_{\mathbf{C}}(X)$ (or simply $\mathcal{A}\operatorname{cocat}(X)$ when there is no danger of confusion), is the $\mathcal{A}$-cocategory of any fibrant-cofibrant replacement of $X$ (see Definitions 1.4 and 1.8).

It will follow immediately from Lemmas 1.10 and 1.11 below that the $\mathcal{A}$-cocategory of an arbitrary (not necessarily fibrant-cofibrant) object is well-defined.

Lemma 1.10. *Let $W$ and $Y$ be two fibrant-cofibrant replacements of an object $X$ in a model category $\mathbf{C}$. Then there is a homotopy equivalence $h\colon W \xrightarrow{\sim} Y$.*

This is a standard result in model category theory, so we omit the proof.

Lemma 1.11. *If $f\colon X \xrightarrow{\sim} Y$ is a weak equivalence of fibrant-cofibrant objects in $\mathbf{C}$, then $\mathcal{A}\operatorname{cocat}(X) = \mathcal{A}\operatorname{cocat}(Y)$.*

Proof. If $\alpha\colon Y \to U \in \mathcal{A}$ is an $\mathcal{A}$-right-versal map for $Y$, then the composite $\beta := \alpha f\colon X \to U^X = U$ is an $\mathcal{A}$-right-versal map for $X$, since $f$ has a homotopy inverse, say, $h$ [**25**, 1.2.8]. Thus one can construct maps $f_n\colon U^X_n \to U_n$ inductively such that the solid-arrow square

$$\begin{array}{ccc} X & \underset{\overset{h}{\longleftarrow \cdots}}{\overset{f}{\longrightarrow}} & Y \\ {\scriptstyle \beta_n}\downarrow & & \downarrow{\scriptstyle \alpha_n} \;\vdots\, {\scriptstyle r} \\ U^X_n & \underset{f_n}{\longrightarrow} & U_n \end{array}$$

is homotopy commutative. If $r\alpha_n \sim \mathrm{id}_Y$, then $hrf_n\beta_n \sim hr\alpha_n f \sim hf \sim \mathrm{id}_X$. So $\mathcal{A}\operatorname{cocat}(X) \leq \mathcal{A}\operatorname{cocat}(Y)$. One can apply the same argument to $h\colon Y \xrightarrow{\sim} X$ and obtain $\mathcal{A}\operatorname{cocat}(Y) \leq \mathcal{A}\operatorname{cocat}(X)$.

This finishes the proof of Lemma 1.11. □

Proposition 1.12. *Suppose $X$ and $Y$ are weakly equivalent objects in $\mathbf{C}$. Then $\mathcal{A}\,\mathrm{cocat}(X) = \mathcal{A}\,\mathrm{cocat}(Y)$. Therefore, $\mathcal{A}$-cocategory is an invariant of the weak equivalence class of an object.*

Proof. It is clearly sufficient to prove the proposition when there is a weak equivalence $f\colon X \to Y$. There is a commutative diagram in $\mathbf{C}$:

$$\begin{array}{ccccc} X & \xrightarrow{r_X} & RX & \xleftarrow{q_{RX}} & QRX \\ {\scriptstyle f}\downarrow & & \downarrow{\scriptstyle Rf} & & \downarrow{\scriptstyle QRf} \\ Y & \xrightarrow[r_Y]{} & RY & \xleftarrow[q_{RY}]{} & QRY \end{array}$$

in which all the horizontal arrows are weak equivalences. Thus, by the 2-out-of-3 axiom of a model category, $QRf$ is a weak equivalence. The result now follows from Lemma 1.11 and Definition 1.9 $\square$

We end this subsection with the following easy but useful observation.

Proposition 1.13. *Let $X$ be an arbitrary object in $\mathbf{C}$. Then $\mathcal{A}\,\mathrm{cocat}\,X = 0$ if and only if its functorial fibrant-cofibrant replacement $QRX$ is a homotopy retract of some object in $\mathcal{A}$.*

**1.3. A convention.** In this short subsection we note that we can assume, without loss of generality, that a right-versal class $\mathcal{A}$ is closed under weak equivalences. If $\mathbf{C}$ is a pointed proper model category and $\mathcal{A}$ is a right-versal class in $\mathbf{C}$ (see Definition 1.3), we consider the two classes $w\mathcal{A}$ and $r\mathcal{A}$, both of which contain $\mathcal{A}$, where:

- $X \in w\mathcal{A}$ if and only if $X$ is fibrant-cofibrant and is weakly equivalent to some object $A \in \mathcal{A}$;
- $X \in r\mathcal{A}$ if and only if $X$ is fibrant-cofibrant and there exist maps $i\colon X \to A$ and $r\colon A \to X$ for some object $A \in \mathcal{A}$ such that the composite $ri$ is homotopic to the identity of $X$.

The following proposition is an immediate consequence of the above definition.

Proposition 1.14. *The classes $r\mathcal{A}$ and $w\mathcal{A}$ are right-versal in $\mathbf{C}$ and $r\mathcal{A}\,\mathrm{cocat}\,X = \mathcal{A}\,\mathrm{cocat}\,X = w\mathcal{A}\,\mathrm{cocat}\,X$ for any object $X$ in $\mathbf{C}$.*

Proposition 1.14 says that the invariant $\mathcal{A}\,\mathrm{cocat}\,X$ will not change if the class $\mathcal{A}$ is replaced by the (possibly bigger) classes $w\mathcal{A}$ or $r\mathcal{A}$. For various reasons, the class $w\mathcal{A}$ is easier to work with than either $\mathcal{A}$ or $r\mathcal{A}$, and so we adopt the following

Convention 1.15. From now on, when we say that a non-empty class $\mathcal{A}$ is a right-versal class in $\mathbf{C}$ we mean:

- every object in $\mathcal{A}$ is fibrant and cofibrant;
- every fibrant-cofibrant object in $\mathbf{C}$ has an $\mathcal{A}$-right-versal map (see Definition 1.3);
- $\mathcal{A}$ is closed under weak equivalences in $\mathbf{C}_{cf}$.

**1.4. Retract property.** The purpose of this subsection is to prove the following retract property of $\mathcal{A}$-cocategory.

Theorem 1.16. *Let $X$ and $Y$ be objects in a pointed proper model category $\mathbf{C}$. If $Y$ is a retract of $X$ (see Definition 1.17 below), then $\mathcal{A}\,\mathrm{cocat}_{\mathbf{C}}(Y) \leq \mathcal{A}\,\mathrm{cocat}_{\mathbf{C}}(X)$.*

Before proving Theorem 1.16 let us make precise what we mean by retracts. Recall that $Q$ and $R$ are, respectively, the functorial cofibrant replacement functor and the functorial fibrant replacement functor in $\mathbf{C}$.

Definition 1.17. Let $X$ and $Y$ be two objects in $\mathbf{C}$. We say that $Y$ is a *retract* of $X$ if there exist maps $i\colon Y \to X$ and $r\colon X \to Y$ such that $QRr \circ QRi \sim \mathrm{id}_{QRY}$.

Theorem 1.16 is an immediate consequence of the following two Lemmas.

Lemma 1.18. *Let $X$ and $Y$ be two fibrant-cofibrant objects in $\mathbf{C}$, and let $i\colon Y \to X$ and $r\colon X \to Y$ be maps such that the composite $ri$ is homotopic to the identity map of $Y$. If $f\colon X \to U$ is an $\mathcal{A}$-right-versal map for $X$, then the map $\alpha := fi\colon Y \to U^Y = U$ is an $\mathcal{A}$-right-versal map for $Y$.*

Lemma 1.19. *With the same hypotheses as in Lemma 1.18, there exists for each integer $k \geq 1$ a map $\varphi_k\colon U_k^Y \to U_k$ such that $\varphi_k \alpha_k \sim f_k i$.*

Assuming Lemmas 1.18 and 1.19 for the moment, let us give the

Proof of Theorem 1.16. There are maps $i\colon QRY \to QRX$ and $r\colon QRX \to QRY$ such that $ri \sim \mathrm{id}_{QRY}$. If $f\colon QRX \to U$ is an $\mathcal{A}$-right-versal map for $QRX$, $\alpha := fi$ as in Lemma 1.18, and if $s\colon U_k \to QRX$ is a homotopy retraction of $f_k$, then it is clear that the composite $rs\varphi_k\colon U_k^{QRY} \to QRY$, with $\varphi_k$ as in Lemma 1.19, is a homotopy retraction of $\alpha_k$. This finishes the proof of Theorem 1.16. □

Now we prove the two Lemmas used in the proof of Theorem 1.16.

Proof of Lemma 1.18. Let $g\colon Y \to A$ be any map with $A \in \mathcal{A}$. We must show that it factors through $\alpha$ up to homotopy. There exists a map $h\colon U \to A$ such that $hf \sim gr$, by the $\mathcal{A}$-right-versality of $f$, and therefore $g \sim gri \sim hfi = h\alpha$. □

Proof of Lemma 1.19. The proof of this Lemma uses the same sort of induction argument and techniques as in the proof of Lemma 1.7, so we leave the details to the interested reader. □

**1.5. A partial ordering on right-versal classes.** What happens if we enlarge the class $\mathcal{A}$? In this subsection we give an answer to this question.

Definition 1.20. Let $\mathbf{C}$ be a pointed proper model category. Define a partial ordering $\leq$ on the class of right-versal classes in $\mathbf{C}$ as follows: If $\mathcal{A}$ and $\mathcal{B}$ are right-versal classes, then we say that $\mathcal{A} \leq \mathcal{B}$ if and only if $\mathcal{A}\,\mathrm{cocat}\,X \leq \mathcal{B}\,\mathrm{cocat}\,X$ for all objects $X$ in $\mathbf{C}$.

The main result of this subsection is the following

THEOREM 1.21. *Let $\mathcal{A}$ and $\mathcal{B}$ be two right-versal classes in* $\mathbf{C}$. *Then the following conditions are equivalent: (1) $\mathcal{A} \leq \mathcal{B}$; (2) $r\mathcal{B} \subseteq r\mathcal{A}$ (see §1.3).*

PROOF. We first prove (1) ⇒ (2). So let $X \in r\mathcal{B}$. Then $X$ is a fibrant-cofibrant object in $\mathbf{C}$ and there exists an object $B \in \mathcal{B}$ of which $X$ is a retract (see Definition 1.17). Since $\mathcal{B}\,\mathrm{cocat}\, B = 0$, it follows from the retract property (see Theorem 1.16) that $\mathcal{A}\,\mathrm{cocat}\, X \leq \mathcal{A}\,\mathrm{cocat}\, B \leq \mathcal{B}\,\mathrm{cocat}\, B = 0$, and so $X$ is a retract of some object in $\mathcal{A}$, i.e., $X \in r\mathcal{A}$. This proves (1) ⇒ (2).

The implication (2) ⇒ (1) follows readily from Lemma 1.22 below. □

LEMMA 1.22. *Let $\mathcal{A}$ and $\mathcal{B}$ be two right-versal classes in* $\mathbf{C}$ *such that $\mathcal{A} \subseteq \mathcal{B}$, let $X$ be a fibrant-cofibrant object in* $\mathbf{C}$, *and let $\alpha\colon X \to U$ and $\beta\colon X \to V$ be, respectively, $\mathcal{A}$-right-versal and $\mathcal{B}$-right-versal maps for $X$. Then for each integer $k \geq 1$ there exists a map $\gamma_k\colon V_k \to U_k$ such that $\gamma_k \beta_k \sim \alpha_k$.*

PROOF. Since this proof involves the same sort of induction argument and techniques as in the proof of Lemma 1.7, we omit the details. □

**1.6. $\mathcal{A}$-category.** The results in the previous subsections are readily dualizable. In this subsection, we give the definitions and corresponding statements for $\mathcal{A}$-category and omit the strictly dual proofs.

DEFINITION 1.23. Let $\mathbf{C}$ be a pointed proper model category, let $\mathcal{A}$ be a non-empty class of fibrant-cofibrant objects in $\mathbf{C}$, and let $X$ be a fibrant-cofibrant object in $\mathbf{C}$. A map $f\colon A \to X$ with $A \in \mathcal{A}$ is called an $\mathcal{A}$-*left-versal map for $X$*, if any map from an object in $\mathcal{A}$ into $X$ factors through $f$ (not necessarily uniquely) up to homotopy. The class $\mathcal{A}$ is said to be *left-versal* in $\mathbf{C}$ if every fibrant-cofibrant object in $\mathbf{C}$ has an $\mathcal{A}$-left-versal map.

For example, given a pointed space $A$, the class $\mathcal{C}(A)$ of $A$-cellular spaces (cf. [**13**, 2.D]) is left-versal, as is the class of all nonempty, possibly infinite wedges of spheres $S^n$ with $n \geq 1$.

As for right-versal classes we adopt the following

CONVENTION 1.24. When we say that a non-empty class $\mathcal{A}$ is a left-versal class in $\mathbf{C}$ we mean: (1) every object in $\mathcal{A}$ is fibrant and cofibrant, (2) every fibrant-cofibrant object in $\mathbf{C}$ has an $\mathcal{A}$-left-versal map, and (3) $\mathcal{A}$ is closed under weak equivalences (in $\mathbf{C}_{cf}$).

Now let $\mathcal{A}$ be a left-versal class in $\mathbf{C}$ and let $f\colon A \to X$ be an $\mathcal{A}$-left-versal map for a fibrant-cofibrant object $X$. We can apply the dual of the construction preceding Definition 1.4 to the map $f$ and obtain a sequence of fibrations $f_n\colon A_n \to X$, $n \geq 1$.

DEFINITION 1.25. Let $\mathcal{A}$ be a left-versal class in $\mathbf{C}$ and let $X$ be a fibrant-cofibrant object in $\mathbf{C}$. The $\mathcal{A}$-*category* of $X$ in $\mathbf{C}$, denoted $\mathcal{A}\,\mathrm{cat}_{\mathbf{C}}(X)$

(or just $\mathcal{A}\,\mathrm{cat}(X)$ if there is no danger of confusion), is the least integer $n \geq 0$ (possibly $\infty$) such that $f_{n+1}$ has a homotopy section, that is, a map $s\colon X \to A_{n+1}$ such that $f_{n+1}s \sim \mathrm{id}_X$. The $\mathcal{A}$-*category* of an arbitrary (not-necessarily fibrant-cofibrant) object is defined to be that of any of its fibrant-cofibrant replacement.

Whenever the symbol $\mathcal{A}\,\mathrm{cat}_{\mathbf{C}}$ (or $\mathcal{A}\,\mathrm{cat}$) appears it is understood that $\mathcal{A}$ is a left-versal class.

As for $\mathcal{A}$-cocategory, the $\mathcal{A}$-category of an arbitrary object $X$ is well-defined.

Proposition 1.26. *Suppose that $X$ and $Y$ are weakly equivalent objects in $\mathbf{C}$. Then $\mathcal{A}\,\mathrm{cat}_{\mathbf{C}}(X) = \mathcal{A}\,\mathrm{cat}_{\mathbf{C}}(Y)$.*

Proposition 1.27. *Let $X$ be an arbitrary object in $\mathbf{C}$. Then $\mathcal{A}\,\mathrm{cat}(X) = 0$ if and only if its functorial fibrant-cofibrant replacement $QRX$ is a homotopy retract of some object in $\mathcal{A}$.*

Theorem 1.28. *Let $Y$ be a retract of $X$ (see Definition 1.17). Then $\mathcal{A}\,\mathrm{cat}(Y) \leq \mathcal{A}\,\mathrm{cat}(X)$.*

Proposition 1.29. *Let $\mathcal{A}$ and $\mathcal{B}$ be two left-versal classes in $\mathbf{C}$. Then $r\mathcal{A} \subseteq r\mathcal{B}$ if and only if $\mathcal{B}\,\mathrm{cat}(W) \leq \mathcal{A}\,\mathrm{cat}(W)$ for all objects $W$ in $\mathbf{C}$.*

## 2. Alternative characterizations of $\mathcal{A}$-(co)category

In the previous section we introduced $\mathcal{A}$-(co)category in a (pointed, proper) model category $\mathbf{C}$. The main result of this section, the Homotopy Pullback (resp. Pushout) Property, gives several alternative characterizations of $\mathcal{A}$-cocategory (resp. $\mathcal{A}$-category). These characterizations will be used repeatedly in §3.

Before stating the main result of this section we need some definitions.

Definition 2.1. Let $\mathcal{A}$ be a right-versal class in $\mathbf{C}$ and let $X$ be a fibrant-cofibrant object in $\mathbf{C}$. The $\mathcal{A}$-*fiber-length* of $X$, denoted $\mathcal{A}\,\mathrm{f.l.}_{\mathbf{C}}(X)$ (or simply $\mathcal{A}\,\mathrm{f.l.}(X)$ if there is no danger of confusion), is defined inductively as follows: $\mathcal{A}\,\mathrm{f.l.}(X) = 0$ if and only if $X$ is in $\mathcal{A}$; $\mathcal{A}\,\mathrm{f.l.}(X) = n$ for some integer $n > 0$ if and only if $\mathcal{A}\,\mathrm{f.l.}(X) \neq m$ for any $0 \leq m < n$ and there is a pullback $X' = \lim\left(Y \xrightarrow{p} Z \xleftarrow{q} W\right)$ in $\mathbf{C}$ such that $q$ is a fibration, $W$ lies in $\mathcal{A}$, $\mathcal{A}\,\mathrm{f.l.}(Y) = n-1$ and $X'$ is weakly equivalent to $X$.

If $X$ is an arbitrary object in $\mathbf{C}$, then the $\mathcal{A}$-*fiber-length* of $X$, denoted $\mathcal{A}\,\mathrm{f.l.}_{\mathbf{C}}(X)$ (or simply $\mathcal{A}\,\mathrm{f.l.}(X)$ if there is no danger of confusion), is defined to be the $\mathcal{A}$-fiber-length of any fibrant-cofibrant replacement of $X$.

It is clear that the $\mathcal{A}$-fiber-length of a fibrant-cofibrant object $X$ is invariant under weak equivalences in $\mathbf{C}_{cf}$, and hence, by Lemma 1.10, the $\mathcal{A}$-fiber-length of an arbitrary object is also invariant under weak equivalences.

From now on, whenever the symbol $\mathcal{A}\,\mathrm{f.\,l.}_{\mathbf{C}}$ (or $\mathcal{A}\,\mathrm{f.\,l.}$) occurs, it is understood that $\mathcal{A}$ is right-versal.

Dually, we can define the $\mathcal{A}$-*cone-length* of $X$, denoted $\mathcal{A}\,\mathrm{c.\,l.}_{\mathbf{C}}(X)$ (or simply $\mathcal{A}\,\mathrm{c.\,l.}(X)$), if $\mathcal{A}$ is left-versal. Whenever the symbol $\mathcal{A}\,\mathrm{c.\,l.}_{\mathbf{C}}$ (or $\mathcal{A}\,\mathrm{c.\,l.}$) occurs, it is understood that $\mathcal{A}$ is left-versal.

Remark 2.2. In [**7, 8**] Cornea studied the cone-length of a space and showed that his cone-length and Ganea's strong category coincide. Our $\mathcal{A}$-cone-length is a generalization of Cornea's: Take $\mathcal{A}$ to be the class of contractible cell-complexes and $\mathcal{A}\,\mathrm{c.\,l.}_{\mathrm{Top}_*}$ is Cornea's cone-length.

We are now ready for the main result of this section:

Theorem 2.3 (Homotopy Pullback Property). *Let $\mathcal{A}$ be a right-versal class in $\mathbf{C}$ and let $X$ be an arbitrary object in $\mathbf{C}$. The following statements are equivalent:*

1. $\mathcal{A}\,\mathrm{cocat}(X) \leq n$.
2. *There exists a fibrant-cofibrant object $Z$ with $\mathcal{A}$-fiber-length $\leq n$ of which $QRX$ is a homotopy retract.*
3. *There exists a fibrant-cofibrant object $Z$ of which $QRX$ is a homotopy retract and such that $Z$ is weakly equivalent to the pullback of a diagram $\left(Z_2 \overset{p}{\twoheadrightarrow} Z_0 \overset{q}{\leftarrow} Z_1\right)$ in which $p$ is a fibration, $\mathcal{A}\,\mathrm{cocat}(Z_1) \leq n-1$, $Z_2$ is weakly equivalent to some object in $\mathcal{A}$.*
4. *Same as (3) with the last condition replaced with $\mathcal{A}\,\mathrm{cocat}(Z_2) = 0$.*
5. *There exist commutative squares*

$$\begin{array}{ccc} X_i & \longrightarrow & Y_i \\ \downarrow & & \downarrow \\ W_i & \longrightarrow & Z_i \end{array} \qquad (0 \leq i \leq n-1)$$

*such that*
   - *$QRX$ is a homotopy retract of $X_{n-1}$,*
   - *$Y_0$ and each $W_i \in \mathcal{A}$,*
   - *$X_j \sim Y_{j+1}$ for $j = 0, \ldots, n-2$,*
   - *each square is a homotopy pullback of fibrant-cofibrant objects.*

Dually, we have

Theorem 2.4 (Homotopy Pushout Property). *Let $\mathcal{A}$ be a left-versal class in $\mathbf{C}$ and let $X$ be an arbitrary object in $\mathbf{C}$. The following statements are equivalent:*

1. $\mathcal{A}\,\mathrm{cat}(X) \leq n$.
2. *There exists a fibrant-cofibrant object $Z$ with $\mathcal{A}$-cone-length $\leq n$ of which $QRX$ is a homotopy retract.*
3. *There exists a fibrant-cofibrant object $Z$ of which $QRX$ is a homotopy retract and such that $Z$ is weakly equivalent to the pushout of a diagram $\left(Z_2 \overset{p}{\leftarrowtail} Z_0 \overset{q}{\to} Z_1\right)$ in which $p$ is a cofibration, $\mathcal{A}\,\mathrm{cat}(Z_1) \leq n-1$, $Z_2$ is weakly equivalent to some object in $\mathcal{A}$.*

4. *Same as (3) with the last condition replaced with* $\mathcal{A}\,\mathrm{cat}(Z_2) = 0$.
5. *There exist commutative squares*

$$\begin{array}{ccc} Z_i & \longrightarrow & W_i \\ \downarrow & & \downarrow \\ Y_i & \longrightarrow & X_i \end{array} \qquad (0 \le i \le n-1)$$

*such that*
- $QRX$ *is a homotopy retract of* $X_{n-1}$,
- $Y_0$ *and each* $W_i \in \mathcal{A}$,
- $X_j \sim Y_{j+1}$ *for* $j = 0, \dots, n-2$,
- *each square is a homotopy pushout of fibrant-cofibrant objects.*

Although the definition of $\mathcal{A}$-cocategory (resp. $\mathcal{A}$-category) involves both homotopy pullback and homotopy pushout squares, the last condition of Theorem 2.3 (resp. 2.4) says that $\mathcal{A}$-cocategory (resp. $\mathcal{A}$-category) can be characterized completely by homotopy pullbacks (resp. pushouts).

Since the proof of Theorem 2.4 is strictly dual to that of Theorem 2.3, we will only give the proof of the later.

Proof of Theorem 2.3. We first prove the equivalence of conditions (1), (3), and (4), for which we need the following three Lemmas.

Lemma 2.5. *Let* $f\colon X \to Y$ *be a map between fibrant-cofibrant objects in* $\mathbf{C}$ *and let* $\alpha\colon X \to U$ *and* $\beta\colon Y \to V$ *be* $\mathcal{A}$*-right-versal maps for* $X$ *and* $Y$*, respectively. Then for each integer* $k \ge 1$ *there exists a map* $f_k\colon U_k \to V_k$ *such that* $f_k\alpha_k \sim \beta_k f$.

Lemma 2.6. *Let*

$$\begin{array}{ccc} Z & \xrightarrow{g} & Z_1 \\ {\scriptstyle f}\downarrow & & \downarrow{\scriptstyle q} \\ Z_2 & \xrightarrow[p]{} & Z_0 \end{array}$$

*be a pullback in* $\mathbf{C}$ *such that* $p$ *is a fibration and* $\mathcal{A}\,\mathrm{cocat}(Z_2) = 0$. *Then* $\mathcal{A}\,\mathrm{cocat}(Z) \le 1 + \mathcal{A}\,\mathrm{cocat}(Z_1)$.

Lemma 2.7. *Let* $f\colon X \to U$ *be an* $\mathcal{A}$*-right-versal map for a fibrant-cofibrant object* $X$ *in* $\mathbf{C}$. *Then* $\mathcal{A}\,\mathrm{cocat}(U_k) < k$ *for each integer* $k \ge 1$.

Assuming Lemmas 2.5, 2.6 and 2.7 for the moment, let us prove the equivalence of conditions (1), (3), and (4).
(1) $\Rightarrow$ (3). Let $f\colon QRX \to U$ be an $\mathcal{A}$-right-versal map for $QRX$. The hypothesis implies that $f_{n+1}$ has a homotopy retraction. But $U_{n+1}$ is weakly equivalent to $U'_{n+1}$, which is the pullback of the diagram $\left(V_n \xrightarrow{\bar{f}_n} U_1 \vee_X U_n \leftarrow U_n\right)$ in which $\bar{f}_n$ is a fibration, $V_n$ is weakly equivalent to $U_1$ and $\mathcal{A}\,\mathrm{cocat}(U_n) \le n-1$ by Lemma 2.7. This proves that (1) $\Rightarrow$ (3).

(3) $\Rightarrow$ (4). This is immediate because every object in $\mathcal{A}$ has $\mathcal{A}$-cocategory 0 and $\mathcal{A}\,\mathrm{cocat}$ is invariant under weak equivalences.
(4) $\Rightarrow$ (1). This follows from the retract property (Theorem 1.16), invariance of $\mathcal{A}$-cocategory under weak equivalences (Proposition 1.12) and Lemma 2.6.

To prove (1) $\Leftrightarrow$ (2), we need the following two Lemmas.

LEMMA 2.8. *Let $f\colon QRX \to U$ be an $\mathcal{A}$-right-versal map for $QRX$. Then $\mathcal{A}\,\mathrm{f.l.}(U_n) < n$ for each integer $n \geq 1$.*

LEMMA 2.9. *For any object $X$ in $\mathbf{C}$, $\mathcal{A}\,\mathrm{cocat}_{\mathbf{C}}(X) \leq \mathcal{A}\,\mathrm{f.l.}_{\mathbf{C}}(X)$.*

Assuming Lemmas 2.8 and 2.9 for the moment, let us prove the equivalence of conditions (1) and (2).
(1) $\Rightarrow$ (2). If $\mathcal{A}\,\mathrm{cocat}\,X \leq n$, then $QRX$ is a homotopy retract of $U_{n+1}$ which by Lemma 2.8 has $\mathcal{A}$-fiber-length at most $n$. This proves (1) $\Rightarrow$ (2).
(2) $\Rightarrow$ (1). There is a pullback $Z' = \lim\left(Y_1 \xrightarrow{p} Y_0 \xleftarrow{q} Y_2\right)$ in which $Z'$ is weakly equivalent to a fibrant-cofibrant object $Z$ of which $QRX$ is a homotopy retract, $q$ is a fibration, $Y_2$ lies in $\mathcal{A}$, and $\mathcal{A}\,\mathrm{f.l.}(Y_1) \leq n-1$. Since $\mathcal{A}\,\mathrm{cocat}(Y_1) \leq n-1$ by Lemma 2.9, it follows from the equivalence of conditions (1) and (3) proved above that $\mathcal{A}\,\mathrm{cocat}(X) \leq (n-1)+1 = n$. This proves (2) $\Rightarrow$ (1).

Now we prove the equivalence of conditions (2) and (5).
(5) $\Rightarrow$ (2). This follows from the easy observation that in condition (5), $\mathcal{A}\,\mathrm{f.l.}(X_i) \leq i+1$ for each $i$.
(2) $\Rightarrow$ (5). There is a pullback $X_{n-1} = \lim\left(W_{n-1} \xrightarrow{q_{n-1}} Z_{n-1} \leftarrow Y_{n-1}\right)$ such that $q_{n-1}$ is a fibration, $W_{n-1} \in \mathcal{A}$, $\mathcal{A}\,\mathrm{f.l.}(Y_{n-1}) \leq n-1$ and $Z \sim X_{n-1}$. Apply the functorial fibrant-cofibrant replacement functor $QR$ to the above square, thus getting a homotopy pullback $QRX_{n-1} = \lim\left(QRW_{n-1} \to QRZ_{n-1} \leftarrow QRY_{n-1}\right)$ of fibrant-cofibrant objects such that $QRX$ is a homotopy retract of $QRX_{n-1}$, $QRW_{n-1} \in \mathcal{A}$, $\mathcal{A}\,\mathrm{f.l.}(QRY_{n-1}) \leq n-1$. Now the same argument applies to $QRY_{n-1}$, $QRY_{n-2}$ etc. and we have shown (2) $\Rightarrow$ (5).

We still have to prove the Lemmas stated above.

PROOF OF LEMMA 2.5. The proof is essentially identical with that of Lemma 1.22, with $\mathrm{id}_X$ replaced with the map $f$. The induction starts because of the $\mathcal{A}$-right-versality of $\alpha$. □

PROOF OF LEMMA 2.6. First we note that, without loss of generality, we may assume that all the four objects in the diagram are fibrant and cofibrant. Indeed, since $Z$ is the (homotopy) limit of the diagram $\{Z_2 \to Z_0 \leftarrow Z_1\}$, there are weak equivalences

$$Z \sim \mathrm{holim}\{QRZ_2 \to QRZ_0 \leftarrow QRZ_1\} \sim Q\,\mathrm{holim}\{QRZ_2 \to QRZ_0 \leftarrow QRZ_1\}.$$

In the last diagram all the objects are fibrant and cofibrant. We can replace $QRp$ by a fibration, all the objects still being fibrant-cofibrant.

Let, then, $\alpha\colon Z \to U$, $\beta\colon Z_1 \to V$, $Z_2 \to T$ and $Z_0 \to S$ be $\mathcal{A}$-right-versal maps for the corresponding objects. Let $r\colon V_k \to Z_1$ be a homotopy retraction of $\beta_k$. We must show that $\alpha_{k+1}$ admits a homotopy retraction. Consider the following diagram in $\operatorname{Ho}\mathbf{C}$:

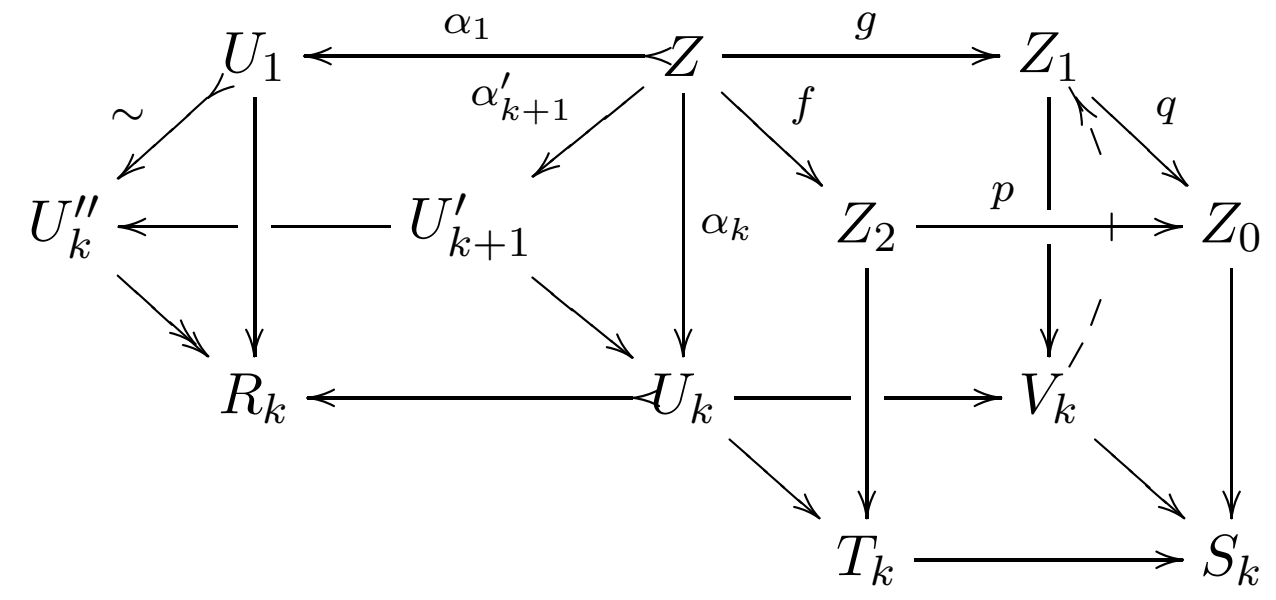

which is commutative, except possibly for the $U_k$-$V_k$-$S_k$-$T_k$ square, and in which $R_k = U_1 \vee_Z U_k$, $U'_{k+1} = U''_k \times_{R_k} U_k$, and $\alpha'_{k+1}$ is the induced map. The rest of the proof consists of diagram-chasing in the above diagram and uses the same techniques as in the proof of Lemma 1.7, so we leave it to the reader. □

Proof of Lemma 2.7. We proceed by induction on $k$, the case $k = 1$ being obvious because $U_1 \in \mathcal{A}$. Suppose that $\mathcal{A}\operatorname{cocat}(U_k) < k$ for some integer $k \geq 1$. Since $\mathcal{A}$-cocategory is invariant under weak equivalences (see Proposition 1.12), it suffices to show that $\mathcal{A}\operatorname{cocat}(U'_{k+1}) \leq k$. There is a pullback $U'_{k+1} = \lim(V_k \twoheadrightarrow R_k \leftarrow U_k)$ in $\mathbf{C}$ in which $V_k$ is weakly equivalent to $U_1$, and hence has $\mathcal{A}$-cocategory 0, and $\mathcal{A}\operatorname{cocat} U_k \leq k-1$ by induction hypothesis. The result now follows from Lemma 2.6. □

Proof of Lemma 2.8. The proof is by induction on $n$. If $n = 1$, then $\mathcal{A}\,\mathrm{f.l.}(U_1) = 0$ because $U_1 \in \mathcal{A}$. Suppose $n > 1$. There is a pullback $U'_n = \lim(U'_1 \twoheadrightarrow U_1 \vee_{QRX} U_{n-1} \leftarrow U_{n-1})$ in which $U'_n$ is weakly equivalent to $U_n$, $U'_1$ is weakly equivalent to $U_1$, and $\mathcal{A}\,\mathrm{f.l.}(U_{n-1}) < n-1$ by induction hypothesis. Thus $\mathcal{A}\,\mathrm{f.l.}(U_n) < n$ by definition. This finishes the induction and the proof of Lemma 2.8. □

Proof of Lemma 2.9. This follows immediately from the equivalence of conditions (1) and (3), whose proof does not use Lemma 2.9. □

The proof of Theorem 2.3 is now complete. □

## 3. Applications of the Homotopy Pullback and Pushout Properties

In this section we apply the main result of §2, the Homotopy Pullback and Pushout Properties, to obtain further properties of our generalized (co)category. The main result of §3.1 is that our generalized (co)category is invariant under Quillen equivalences which are also modelization functors.

In §3.2 we show that our general approach to (co)category can be specialized to obtain a reasonable notion of "rational cocategory" of spaces which coincides with Ganea's cocategory for (nice) rational spaces and which also satisfies the fibration property. The final section contains formulae for the generalized cocategory of a homotopy function complex and the generalized cocategory (resp. category) of a finite product (resp. coproduct).

**3.1. $\mathcal{A}$-(co)category under Quillen functors.** In this subsection we show that $\mathcal{A}$-(co)category is not changed under Quillen modelization equivalences, i.e. Quillen equivalences (see [**25**]) which are also modelization functors (see [**10**]). Recall that a modelization functor is a functor that preserves weak equivalences, homotopy pullbacks and homotopy pushouts.

We begin with a preliminary result.

Proposition 3.1. *Let $F\colon \mathbf{C} \to \mathbf{D}$ be a (left or right) Quillen functor that is also a modelization functor. Let $\mathcal{A}$ be a right-versal class in $\mathbf{C}$ and $\mathcal{B}$ be a right-versal class in $\mathbf{D}$ such that, whenever $Y \in \mathbf{D}_{cf}$ is weakly equivalent to $FA$ for some object $A \in \mathcal{A}$, then $Y \in \mathcal{B}$. Then for any object $X$ in $\mathbf{C}$, $\mathcal{A}\,\mathrm{cocat}_{\mathbf{C}}(X) \geq \mathcal{B}\,\mathrm{cocat}_{\mathbf{D}}(FX)$.*

Now we can state the result concerning the invariance of $\mathcal{A}$-cocategory under Quillen modelization equivalences.

Theorem 3.2. *Let $F\colon \mathbf{C} \rightleftarrows \mathbf{D}\colon U$ be a pair of Quillen equivalences with $F$ and $U$ both modelization functors. Let $\mathcal{A}$ be a right-versal class in $\mathbf{C}$ and $\mathcal{B}$ be a right-versal class in $\mathbf{D}$. Assume that $F$ takes fibrant objects to fibrant objects and $U$ takes cofibrant objects to cofibrant objects. Then:*

1. *There exists a right-versal class $\mathcal{A}'$ in $\mathbf{D}$ such that for any object $X$ in $\mathbf{C}$, $\mathcal{A}\,\mathrm{cocat}_{\mathbf{C}}(X) = \mathcal{A}'\,\mathrm{cocat}_{\mathbf{D}}(FX)$.*
2. *There exists a right-versal class $\mathcal{B}'$ in $\mathbf{C}$ such that for any object $Y$ in $\mathbf{D}$, $\mathcal{B}\,\mathrm{cocat}_{\mathbf{D}}(Y) = \mathcal{B}'\,\mathrm{cocat}_{\mathbf{C}}(UY)$.*

This theorem has an Eckmann-Hilton dual in which right-versal is replaced by left-versal and $\mathcal{A}\,\mathrm{cocat}$ is replaced by $\mathcal{A}\,\mathrm{cat}$ throughout. Its proof is strictly dual to that of 3.2. We will not state this result explicitly.

Theorem 3.2 and its dual apply, for example, to the adjoint pair

$$|-|\colon \mathrm{SS}_* \rightleftarrows \mathrm{Top}_*\colon \mathrm{Sing}$$

of Quillen equivalences, where $|-|$ and Sing denote, respectively, the geometric realization functor and the singular complex functor.

Proof of Proposition 3.1. Write $n = \mathcal{A}\,\mathrm{cocat}_{\mathbf{C}}\,X$. If $n = \infty$ then there is nothing to prove, and the case $n = 0$ is easy. So suppose $0 < n < \infty$. Then there are homotopy pullback diagrams as in Theorem 2.3. Applying the modelization functor $QRF$ to them we see at once that $\mathcal{B}\,\mathrm{cocat}_{\mathbf{D}}(FX_{n-1}) = \mathcal{B}\,\mathrm{cocat}_{\mathbf{D}}(QRFX_{n-1}) \leq n$. Since $FQRX$ is a retract of $FX_{n-1}$, $\mathcal{B}\,\mathrm{cocat}_{\mathbf{D}}(FX) = \mathcal{B}\,\mathrm{cocat}_{\mathbf{D}}(FQRX) \leq n$. This finishes the proof of Proposition 3.1. □

Proof of Theorem 3.2. We need the following two Lemmas.

Lemma 3.3. *Let $F\colon \mathbf{C} \rightleftarrows \mathbf{D}\colon U$ be a pair of Quillen equivalences and let $\mathcal{A}$ and $\mathcal{B}$ be right-versal classes in $\mathbf{C}$ and $\mathbf{D}$, respectively. Assume that $F$ takes fibrant objects to fibrant objects and $U$ takes cofibrant objects to cofibrant objects. Then:*

1. *The class $\mathcal{A}'$ consisting of all fibrant-cofibrant objects $Y$ in $\mathbf{D}$ such that there exists an object $A \in \mathcal{A}$ with $Y$ weakly equivalent to $FA$, is right-versal in $\mathbf{D}$.*
2. *The class $\mathcal{B}'$ consisting of all fibrant-cofibrant objects $X$ in $\mathbf{C}$ such that there exists an object $B \in \mathcal{B}$ with $X$ weakly equivalent to $UB$, is right-versal in $\mathbf{C}$.*

Lemma 3.4. *With the same notations and assumptions as in Lemma 3.3, let*

- *$\mathcal{A}''$ be the class consisting of all fibrant-cofibrant objects $Z$ in $\mathbf{C}$ such that there exists an object $M \in \mathcal{A}'$ with $Z$ weakly equivalent to $UM$, and*
- *$\mathcal{B}''$ be the class consisting of all fibrant-cofibrant objects $W$ in $\mathbf{D}$ such that there exists an object $N \in \mathcal{B}'$ with $W$ weakly equivalent to $FN$.*

*Then $\mathcal{A} = \mathcal{A}''$ and $\mathcal{B} = \mathcal{B}''$.*

Assuming Lemmas 3.3 and 3.4, let us finish the proof of Theorem 3.2. Since the proofs of the two parts are very similar, we only give the proof of part 1. Let $\mathcal{A}'$ and $\mathcal{A}''$ be as in Lemmas 3.3 and 3.4. We then have

$$\mathcal{A}\,\mathrm{cocat}_{\mathbf{C}}(X) \overset{(a)}{\geq} \mathcal{A}'\,\mathrm{cocat}_{\mathbf{D}}(FX) \overset{(b)}{\geq} \mathcal{A}''\,\mathrm{cocat}_{\mathbf{C}}(UFX)$$

$$\overset{(c)}{=} \mathcal{A}\,\mathrm{cocat}_{\mathbf{C}}(UFX) \overset{(d)}{=} \mathcal{A}\,\mathrm{cocat}_{\mathbf{C}}(X).$$

The inequalities (a) and (b) follow from Proposition 3.1, (c) follows from Lemma 3.4, and (d) follows from the weak equivalences

$$X \sim \mathbf{R}U\mathbf{L}FX = URFQX \sim UFX.$$

Modulo the proofs of Lemmas 3.3 and 3.4, the proof of Theorem 3.2 is complete. □

Now we give the proofs of Lemmas 3.3 and 3.4.

Proof of Lemma 3.3. We will prove the first part; the proof of the second part is similar.

Let $Y$ be an object of $\mathbf{D}$ that is both fibrant and cofibrant. We must show that it has an $\mathcal{A}'$-right-versal map. By using the fact that $\mathbf{L}F\colon\ \mathrm{Ho}\,\mathbf{C} \to \mathrm{Ho}\,\mathbf{D}$ is an equivalence of categories, it is not hard to see that there exists a weak equivalence $\alpha\colon Y \xrightarrow{\sim} FZ$ in $\mathbf{D}$ with $Z$ fibrant and cofibrant in $\mathbf{C}$.

The following claim will finish the proof of the Lemma.

*Claim 1.* The composite $\gamma := F\beta \circ \alpha\colon Y \to FA$, where $\beta\colon Z \to A$ is an $\mathcal{A}$-right-versal map for $Z$, is an $\mathcal{A}'$-right-versal map for $Y$.

So let $\tau\colon Y \to W$ be a map in $\mathbf{D}$ with $W$ fibrant-cofibrant and is weakly equivalent to $FB$ for some object $B$ in $\mathcal{A}$. We must show that $\tau$ factors through $\gamma$ up to homotopy. There is a weak equivalence $w\colon FB \to W$ in $\mathbf{D}$, since both $W$ and $FB$ are fibrant and cofibrant. To finish the proof of Claim 1, we will construct a map $\varphi\colon A \to B$ in $\mathbf{C}$ with the property that $w \circ F\varphi \circ \gamma \sim \tau$. It follows from the hypotheses that the unit maps $Z \to \mathbf{R}U\mathbf{L}FZ$ and $B \to \mathbf{R}U\mathbf{L}FB$ are isomorphisms in $\operatorname{Ho}\mathbf{C}$, and it follows from Ken Brown's Lemma (cf. [**25**, 1.1.12]) that the natural maps $\varepsilon_Z\colon Z \to UFZ$ and $\varepsilon_B\colon B \to UFB$ are weak equivalences and hence homotopy equivalences (cf. [**25**, 1.2.8]). By the $\mathcal{A}$-right-versality of $\beta$, there exists a map $\varphi\colon A \to B$ in $\mathbf{C}$ such that $\varphi\circ\beta \sim \psi$, where $\psi$ is the composite $\varepsilon_B^{-1}\circ U(w^{-1}\tau\alpha^{-1})\circ\varepsilon_Z\colon Z \to B$. It is now easy to see that the map $\varphi$ has the desired property.

The proof of Claim 1, and hence of Lemma 3.3, is complete. $\square$

Proof of Lemma 3.4. We will prove part 1; the proof of part 2 is similar.

First let $X$ be an object in $A''$. Then $X$ is a fibrant-cofibrant object in $\mathbf{C}$ and is weakly equivalent to $UW$ for some object $W$ in $\mathcal{A}'$. So $W$ is fibrant-cofibrant in $\mathbf{D}$ and is weakly equivalent to $FA$ for some object $A$ in $\mathcal{A}$. Thus there are weak equivalences $f\colon X \xrightarrow{\sim} UW$ in $\mathbf{C}$ and $g\colon W \xrightarrow{\sim} FA$ in $\mathbf{D}$ (cf. [**25**, 1.2.8]). By Ken Brown's Lemma the map $Ug\colon UW \to UFA$ is a weak equivalence in $\mathbf{C}$, and therefore so is the composite $Ug\circ f\colon X \to UFA$. Also, it follows from our assumption and Ken Brown's Lemma that there are weak equivalences

$$A \longrightarrow URFQA \xleftarrow{Ur_{FQA}} UFQA \xrightarrow{UFq_A} UFA,$$

and so $X$ is weakly equivalent to $A$. Hence $X$ is in $\mathcal{A}$.

Conversely, let $B$ be in $\mathcal{A}$. Then, as above, $B$ is weakly equivalent to $UFB$. But it follows again from our assumption that $FB$ is in $A'$, and so $B$ is in $A''$. This finishes the proof of Lemma 3.4. $\square$

**3.2. Rational LS cocategory.** Unless otherwise specified, all spaces (and simplicial sets) in this subsection are pointed, simply connected, and have finite $\mathbf{Q}$-type (i.e. finite dimensional rational homology group in each dimension).

Before stating the results of this subsection let us briefly review the subject of rational cocategory. The work of Felix and Halperin [**14**] on rational category (i.e. $\mathrm{cat}_0$) inspired Sbaï [**28, 29**] to introduce rational cocategory in terms of the Quillen minimal model of a space, by using (the dual of) one of the characterizations of rational category. He showed that his rational cocategory is an upper-bound for Ganea's cocategory for a rational space, and conjectured that equality holds in general. This conjecture was disproved by M. Hovey [**24**]. The more recent work of Doeraene [**10**], as far as cocategory is concerned, applies only to what he calls $J^{op}$-categories; the category $\mathrm{Top}_*$ of pointed spaces, however, is not a $J^{op}$-category. This leaves open the question of what the right definition of rational cocategory should be.

In this subsection we give an alternative definition of rational cocategory of a space, which we denote by $\mathrm{cocat}_0(X)$, in terms of its Sullivan minimal model, and show that our rational cocategory has many of the properties that one would like it to have.

Here is our definition of rational cocategory:

DEFINITION 3.5. The *rational cocategory* of a space $X$, denoted $\mathrm{cocat}_0(X)$, is defined to be $\mathrm{cat}_{\mathrm{DGA}_*}(AX)$.

Hereafter whenever we use the term "rational cocategory" we are referring to Definition 3.5.

In the definition above, $A$ is the Sullivan-de Rham functor from pointed simplicial sets $\mathrm{SS}_*$ to $\mathrm{DGA}_*$, the category of augmented commutative differential graded algebras over the rationals with Bousfield-Gugenheim's model category structure [**2**], and $\mathrm{(co)cat}_{\mathrm{DGA}_*}$ is the (co)category in $\mathrm{DGA}_*$ as defined in §1 (Definitions 1.9 and 1.25 with $\mathcal{A} = \{*\}$).

The motivation for the above definition is the following observation.

PROPOSITION 3.6. *The equality* $\mathrm{cat}_0(X) = \mathrm{cocat}_{\mathrm{DGA}_*}(AX)$ *holds for any space* $X$.

REMARK 3.7. In [**10**, p. 259] Doeraene made a similar remark that $\mathrm{cat}_0(X) = {}_{\mathrm{CDA}^{*c0}}\,\mathrm{cocat}(X)$, where the right-hand side is cocategory (in the sense of Doeraene) in the $J^{op}$-category $\mathrm{CDA}^{*c0}$ of augmented commutative $c$-connected differential graded algebras. However, $\mathrm{CDA}^{*c0}$ is not a $J$-category (cf. [**10**, p. 257]), and so the results in [**10**] do not immediately apply to LS category in $\mathrm{CDA}^{*c0}$ (or $\mathrm{DGA}_*$).

Our main result here states that Ganea's cocategory rationalizes to our rational cocategory. Thus, the latter does not have the disadvantage of Sbaï's rational cocategory.

THEOREM 3.8. *Let* $X$ *be a space. Then* $\mathrm{cocat}(X) \geq \mathrm{cocat}_0(X)$. *Moreover, equality holds if* $X$ *is also a rational space. In particular,* $\mathrm{cocat}(X_{\mathbf{Q}}) = \mathrm{cocat}_0(X)$.

Here $X_{\mathbf{Q}}$ denotes the rationalization of $X$ and cocat means Ganea's cocategory. Since the cocategories of a topological space and of its singular complex coincide (see Theorem 3.2), there is no ambiguity in leaving out the subscript $\mathrm{SS}_*$ in cocat.

As mentioned above, Hovey disproved the conjecture that Ganea's cocategory rationalizes to Sbaï's rational cocategory. He actually did this by exhibiting a fibration $(F \to E \to B)$ of rational spaces in which $E$ is an $H$-space, hence has cocategory 1, while $F$ has Sbaï's rational cocategory $> 2$. Since cocategory in the sense of Ganea satisfies the fibration property–that in a fibration the fiber has cocategory no bigger than that of the total space plus 1–this fibration disproved Sbaï's conjecture. The next result shows that our rational cocategory does satisfy the fibration property (and more).

PROPOSITION 3.9 (Rational fibration property). *Let* $\mathcal{E} = (F \to E \to B)$ *be a rational fibration in the sense of Halperin* [**20**] *with* $B$ *connected. Then* $\mathrm{cocat}_0(F) \leq 1 + \mathrm{cocat}_0(E)$. *In particular, this inequality holds if* $\mathcal{E}$ *is a fibration.*

REMARK 3.10. Recall that a rationally coformal space $X$ is a rational space such that $\pi_*(\Omega X) \otimes \mathbf{Q}$ equipped with zero differential has the same Quillen minimal model as $X$. Sbaï [**28, 29**] showed that when $X$ is rationally coformal, his rational cocategory of $X$ coincides with $\mathrm{cocat}(X)$, and hence, by Theorem 3.8, also with our $\mathrm{cocat}_0(X)$.

Now we give the proofs of the results above.

PROOF OF PROPOSITION 3.6. Recall from [**14**] that the *rational category* of $X$ is the least integer $n \geq 0$ (possibly $\infty$) such that the Sullivan minimal model $\Lambda V$ of $X$ is a homotopy retract of the rational Ganea space $\Gamma_n(\Lambda V)$, which has the rational homotopy type of the space $B_n\Omega X$ in the Milnor filtration of $B\Omega X$. The proposition is now an immediate consequence of the fact that the functor $M\colon \mathrm{SS}_* \to \mathrm{DGA}_*$, defined as $A$ followed by functorial cofibrant replacement, takes rational fibrations in $\mathrm{SS}_*$ to cofibrations in $\mathrm{DGA}_*$. □

PROOF OF THEOREM 3.8. First we note that if two spaces $X$ and $Y$ have the same rational homotopy type, then $\mathrm{cocat}_0(X) = \mathrm{cocat}_0(Y)$ (see Proposition 1.26).

Now we prove the first assertion of Theorem 3.8.

LEMMA 3.11. *Let* $X$ *be a space. Then* $\mathrm{cocat}(X) \geq \mathrm{cocat}_0(X)$.

PROOF. Consider Ganea's fiber-cofiber construction (GFC) for $X$:

$$\begin{array}{ccccc} & & \vdots & & \\ & & \downarrow{\scriptstyle i_1} & & \\ & & F_1X & \xrightarrow{g_1} & B_1X \\ & \overset{f_1}{\nearrow} & \downarrow{\scriptstyle i_0} & & \\ X & \xrightarrow[f_0]{} & F_0X & \xrightarrow[g_0]{} & B_0X \end{array} \tag{GFC}$$

in which $(X \to F_0X \to B_0X) = (X \to CX \to \Sigma X)$, each sequence $\mathcal{C}_n\colon (X \to F_nX \to B_nX)$ is a cofibration sequence, and each sequence $\mathcal{E}_n\colon (F_{n+1}X \to F_nX \to B_nX)$ is a fibration sequence. Thus, in particular, $\mathrm{cocat}(X)$ is the smallest integer $n$ for which $f_n$ has a homotopy retraction.

We observe the following:

*Claim 1.* Every space in the diagram (GFC) is simply-connected.

*Claim 2.* For $n \geq 0$, $(F_nX)_{\mathbf{Q}} \simeq F_nX_{\mathbf{Q}}$ and $(B_nX)_{\mathbf{Q}} \simeq B_nX_{\mathbf{Q}}$.

*Claim 3.* Every space in the diagram (GFC) has finite $\mathbf{Q}$-type.

*Claim 4.* $M(\mathrm{GFC})$ is precisely the diagram that defines $\mathrm{cat}_{\mathrm{DGA}_*}(AX)$.

Assuming these claims for the moment, notice that a homotopy retraction of $f_n$ yields a homotopy section of $Mf_n$, which by Claim 4 implies $\mathrm{cocat}_0(X) = \mathrm{cat}_{\mathrm{DGA}_*}(AX) \le n$, completing the proof of Lemma 3.11.

To prove Claim 1, note that in any cofibration sequence $A \to X \to X/A$ the homotopy fiber $\mathrm{Fib}(X \to X/A)$ is $A$-cellular (cf. [**4**, 5.4(2)]). Thus each $F_nX$ is simply-connected, since $F_nX = \mathrm{Fib}(F_{n-1}X \to (F_{n-1}X)/X)$. A more elementary way to see this, as the referee suggests, is to use the Blakers-Massey Theorem. That the $B_nX$ are all simply-connected follows from the long exact homotopy sequence of $\mathcal{E}_n$. This finishes the proof of Claim 1.

Claim 2 follows from Claim 1 and ([**31**, §2.9, Prop. 2.4] or [**21**, Ch. II, Cor. 1.10 and 1.11]), which says that for simply-connected spaces rationalization preserves fibration and cofibration sequences.

Claim 3 is proved by induction on $n$, the case $n = 0$ being trivial since $F_0X = CX \simeq *$ and $B_0X = \Sigma X$. Suppose $F_{n-1}X$ and $B_{n-1}X$ have finite $\mathbf{Q}$-type. By Claim 2 we only need establish that $F_nX_{\mathbf{Q}}$ and $B_nX_{\mathbf{Q}}$ have finite type integral homology. For $F_nX_{\mathbf{Q}}$ this follows from induction hypothesis, the long exact homotopy sequence of the fibration $(\mathcal{E}_{n-1})_{\mathbf{Q}}$ and Serre's $\mathcal{C}$-theory [**30**]. The case for $B_nX_{\mathbf{Q}}$ now follows from the long exact homology sequence of the cofibration sequence $(\mathcal{C}_n)_{\mathbf{Q}}$. This finishes the induction and the proof of Claim 3.

Finally, Claim 4 is true since $M(F_0X) \simeq \mathbf{Q}$, $M\mathcal{C}_n$ is a fibration sequence (cf. [**2**, p. 82]), and $M\mathcal{E}_n$ is a cofibration sequence because by Claims 1 and 3 and [**19**, 20.3], $\mathcal{E}_n$ is a rational fibration.

This finishes the proof of Lemma 3.11. □

The second assertion of Theorem 3.8 follows from the following Lemma.

Lemma 3.12. *Let $X$ be a rational space. Then* $\mathrm{cocat}(X) \le \mathrm{cocat}_0(X)$.

Proof. We may assume without loss of generality that $X$ is a Kan complex, since the map $X \xrightarrow{\sim} \mathrm{Sing}\,|X|$ induces a weak equivalence $A(\mathrm{Sing}\,|X|) \xrightarrow{\sim} A(X)$. We claim that for $Y \in \mathrm{DGA}_*$, one has the inequality

$$\mathrm{cocat}(FY) \le \mathrm{cat}_{\mathrm{DGA}_*}(Y),$$

where $F\colon \mathrm{DGA}_* \to \mathrm{SS}_*$ is the functor defined by Bousfield and Gugenheim [**2**]. To see this, observe that the functor $F$ takes weak equivalences between cofibrant algebras to weak equivalences between Kan complexes and takes pushouts to pullbacks. Therefore, the claim follows from the Homotopy Pullback and Homotopy Pushout Properties (cf. Theorems 2.3 and 2.4).

Now since $X$ is a rational space, there is a weak-equivalence $X \xrightarrow{\sim} FMX$ (cf. [**2**, 10.1]). Thus, with $Y = MX$ in the inequality above it follows that $\mathrm{cocat}(X) = \mathrm{cocat}(FMX) \le \mathrm{cat}_{\mathrm{DGA}_*}(MX) = \mathrm{cocat}_0(X)$.

This finishes the proof of Lemma 3.12. □

The proof of Theorem 3.8 is now complete. □

Proof of Proposition 3.9. Recall from [**20**] that a sequence $\mathcal{E}\colon F \xrightarrow{i} E \xrightarrow{\pi} B$ of path-connected spaces is a *rational fibration* if the map $\pi \circ i$ sends $F$ to the base point of $B$ and in the diagram $(*)$ below,

$$\begin{array}{ccccc} A(B) & \xrightarrow{A(\pi)} & A(E) & \xrightarrow{A(i)} & A(F) \\ \| & & \sim\uparrow & & \uparrow\alpha \\ A(B) & \xrightarrow[\varphi]{} & A(B)\otimes \Lambda V & \longrightarrow & \Lambda V \end{array}$$

where $\varphi$ is the minimal model of $A(\pi)$, the induced map $\alpha$ is a weak equivalence.

Now by [**1**, I.8.18] the bottom row of the diagram $(*)$ is a cofibration in $\mathrm{DGA}_*$. Thus it follows that

$$\begin{aligned} \mathrm{cocat}_0(F) &= \mathrm{cat}_{\mathrm{DGA}_*}(A(F)) = \mathrm{cat}_{\mathrm{DGA}_*}(\Lambda V) \\ &\leq 1 + \mathrm{cat}_{\mathrm{DGA}_*}(A(B)\otimes \Lambda V) = 1 + \mathrm{cocat}_0(E). \end{aligned}$$

This proves the first assertion of Theorem 3.9.

The second assertion follows from [**18**] or [**19**, 20.6], which says that under the stated hypotheses $\mathcal{E}$ is a rational fibration.

This completes the proof of Proposition 3.9. □

**3.3. Further applications.** In this final subsection we give two more results concerning our generalized (co)category, both of which are applications of the Homotopy Pullback and Pushout Properties. Our first result here is a generalization of the corresponding classical result of Ganea [**17**].

Theorem 3.13. *Let* $\mathbf{C}$ *be a pointed proper model category in which the weak equivalences are closed under finite products, and let* $\mathcal{A}$ *be a right-versal class in* $\mathbf{C}$ *that is closed under finite products. Then for any objects* $X_1, \ldots, X_n$ *in* $\mathbf{C}$ *we have*

$$\mathcal{A}\,\mathrm{cocat}_{\mathbf{C}}(X_1 \times \cdots \times X_n) = \max\left(\mathcal{A}\,\mathrm{cocat}_{\mathbf{C}}(X_i)\colon 1 \leq i \leq n\right).$$

This theorem has an Eckmann-Hilton dual in which products, right-versal, and $\mathcal{A}\,\mathrm{cocat}$ are replaced by, respectively, coproducts, left-versal, and $\mathcal{A}\,\mathrm{cat}$ throughout. Its proof is strictly dual to that of Theorem 3.13. We will not state it explicitly.

Proof of Theorem 3.13. We will omit the subscript $\mathbf{C}$ in this proof, since $\mathbf{C}$ is the only model category under consideration. It is clearly sufficient to prove the theorem for $n = 2$; the general case follows by an induction argument. So let $X$ and $Y$ be objects in $\mathbf{C}$ and write $m = \mathcal{A}\,\mathrm{cocat}(X)$, $n = \mathcal{A}\,\mathrm{cocat}(Y)$. Since both $X$ and $Y$ are retracts of $X \times Y$, the retract property (Theorem 1.16) shows that both $\mathcal{A}\,\mathrm{cocat}(X)$ and $\mathcal{A}\,\mathrm{cocat}(Y)$ are less than or equal to $\mathcal{A}\,\mathrm{cocat}(X \times Y)$.

It remains to show that $\mathcal{A}\,\mathrm{cocat}(X \times Y)$ is no bigger than the maximum of $\mathcal{A}\,\mathrm{cocat}(X)$ and $\mathcal{A}\,\mathrm{cocat}(Y)$. If either $m$ or $n$ is $\infty$, then there is nothing to prove. So we may assume that both $m$ and $n$ are finite. Note that to obtain the desired inequality we may replace $X$ and $Y$ by, respectively, $QRX$ and $QRY$. In the case that either $m$ or $n$ is zero (or both), the proof is very easy and we omit it. So we suppose that both $m$ and $n$ are positive and finite. There exists an object $Z$ (resp. $W$) of which $QRX$ (resp. $QRY$) is a homotopy retract and such that $Z$ (resp. $W$) is weakly equivalent to the pullback of a diagram $\left(Z_2 \xrightarrow{p} Z_0 \xleftarrow{q} Z_1\right)$ $\left(\text{resp. } \left(W_2 \xrightarrow{r} W_0 \xleftarrow{s} W_1\right)\right)$ such that $p$ (resp. $r$) is a fibration, $\mathcal{A}\,\mathrm{cocat}(Z_1) = m-1$ (resp. $\mathcal{A}\,\mathrm{cocat}(W_1) = n-1$) and $\mathcal{A}\,\mathrm{cocat}(Z_2) = 0 = \mathcal{A}\,\mathrm{cocat}(W_2)$. Thus $QRX \times QRY$ is a homotopy retract of $Z \times W$, which by our hypothesis is weakly equivalent to $\lim\left(Z_2 \times W_2 \xrightarrow{p\times r} Z_0 \times W_0 \xleftarrow{q\times s} Z_1 \times W_1\right)$ in which $\mathcal{A}\,\mathrm{cocat}(Z_2 \times W_2) = 0$, $p \times r$ is a fibration, and $\mathcal{A}\,\mathrm{cocat}(Z_1 \times W_1) \leq \max\{m-1, n-1\}$ by induction hypothesis (on $m+n$). Thus, it follows from the Homotopy Pullback Property (Theorem 2.3) that

$$\mathcal{A}\,\mathrm{cocat}(QRX \times QRY) \leq 1 + \max\{m-1, n-1\} = \max\{m, n\}.$$

This completes the proof of Theorem 3.13. □

Our last application of the Homotopy Pullback and Pushout Properties is the following result. We use the notation $\mathrm{map}_*(-,-)$ for (pointed) homotopy function complex (see [**22**]). Recall that $\mathrm{SS}_*$ denotes the category of pointed simplicial sets.

Theorem 3.14. *Let* $\mathbf{C}$ *be a simplicial pointed proper model category, let* $\mathcal{A}$ *and* $\mathcal{B}$ *be, respectively, left and right-versal classes in* $\mathbf{C}$ *and let* $X$ *be a cofibrant object and* $Y$ *be a fibrant object in* $\mathbf{C}$. *Let* $\mathcal{C}$ *be a right-versal class in* $\mathrm{SS}_*$ *such that:*

- $\mathrm{map}_*(A, Y) \in \mathcal{C}$ *for every object* $A$ *that is weakly equivalent to some object in* $\mathcal{A}$*;*
- $\mathrm{map}_*(X, B) \in \mathcal{C}$ *for every object* $B$ *that is weakly equivalent to some object in* $\mathcal{B}$.

*Then we have*

$$\mathcal{C}\,\mathrm{cocat}_{\mathrm{SS}_*}\,\mathrm{map}_*(X, Y) \leq \min\left(\mathcal{A}\,\mathrm{cat}_{\mathbf{C}}(X), \mathcal{B}\,\mathrm{cocat}_{\mathbf{C}}(Y)\right).$$

Proof. We first consider the inequality involving $\mathcal{A}\,\mathrm{cat}(X)$ (the subscript $\mathbf{C}$ is omitted). Write $m = \mathcal{A}\,\mathrm{cat}(X)$. The case $m = 0$ follows from the assumption and the case $m = \infty$ is trivial. So we now do induction on $m \geq 0$; the case $m = 0$ is already verified. Suppose $m > 0$. By the Homotopy Pushout Property, Theorem 2.4, there exists a fibrant-cofibrant object $Z$ of which $QRX$ is a homotopy retract and $Z$ is weakly equivalent to the homotopy pushout $Z'$ of the diagram $\left(Z_2 \xleftarrow{p} Z_0 \xrightarrow{q} Z_1\right)$ in which $p$ is a cofibration, $Z_2$ is weakly equivalent to some object in $\mathcal{A}$, and $\mathcal{A}\,\mathrm{cat}(Z_1) = m-1$. Applying the functor $\mathrm{map}_*(-, Y)$ to the above homotopy pushout square, it is not difficult to see that the following homotopy pullback diagram in $\mathrm{SS}_*$

(cf. [**11**, 62.1]),

$$\begin{array}{ccc} \operatorname{map}_*(Z', Y) & \longrightarrow & \operatorname{map}_*(Z_1, Y) \\ \downarrow & & \downarrow \\ \operatorname{map}_*(Z_2, Y) & \twoheadrightarrow & \operatorname{map}_*(Z_0, Y), \end{array}$$

together with Theorem 2.3 yield the desired inequality. The proof of the inequality involving $\mathcal{B}\operatorname{cocat}_{\mathbf{C}}(Y)$ is very similar to the argument above, so we omit the details. □

## References

[1] H. Joachim Baues, *Algebraic homotopy*, Cambridge University Press, Cambridge, 1988.

[2] A. K. Bousfield and V. K. A. M. Gugenheim, *On PL de Rham theory and rational homotopy type*, Mem. Amer. Math. Soc. **179**, 1976.

[3] A. K. Bousfield and D. M. Kan, *Homotopy limits, completions and localizations*, Lecture Notes in Math. **304**, Springer-Verlag, 1972.

[4] Wojciech Chachólski, *A generalization of the triad theorem of Blakers-Massey*, Topology **36** (1997) 1381-1400.

[5] Monica Clapp and Dieter Puppe, *Invariants of the Lusternik-Schnirelmann type and the topology of critical sets*, Trans. Amer. Math. Soc. **298** (1986) 603-620.

[6] ———, *The generalized Lusternik-Schnirelmann category of a product space*, Trans. Amer. Math. Soc. **321** (1990) 525-532.

[7] Octavian Cornea, *Cone-length and Lusternik-Schnirelmann category*, Topology **33** (1994) 95-111.

[8] ———, *Strong LS category equals cone-length*, Topology **34** (1995) 377-381.

[9] J.-P. Doeraene, *LS-catégorie dans une catégorie à modèles*, Thesis, Université Catholique de Louvain, 1990.

[10] ———, *L.S.-category in a model category*, J. Pure Appl. Alg. **84** (1993) 215-261.

[11] W. G. Dwyer, P. S. Hirschhorn, and D. M. Kan, *Model categories and more general abstract homotopy theory*, in preparation, available at http://www-math.mit.edu/~psh.

[12] W. G. Dwyer and J. Spalinski, *Homotopy theories and model categories*, in: Handbook of algebraic topology, ed. I. M. James, Elsevier Science, 1995, 73-126.

[13] E. D. Farjoun, *Cellular spaces, null spaces and homotopy localization*, Lecture Notes in Math. **1622**, Springer-Verlag, 1996.

[14] Y. Felix and S. Halperin, *Rational L.S. category and its applications*, Trans. Amer. Math. Soc. **273** (1982) 1-37.

[15] Tudor Ganea, *A generalization of the homology and homotopy suspension*, Comment. Math. Helv. **39** (1965) 295-322.

[16] ———, *Fibrations and cocategory*, Comment. Math. Helv. **35** (1961) 15-24.

[17] ———, *Lusternik-Schnirelmann category and cocategory*, Proc. London Math. Soc. (3) **10** (1960) 623-639.

[18] P. Grivel, Suite spectrale et modele minimal d'une fibration, Thèse, Université de Genève, 1977.

[19] Stephen Halperin, *Lectures on minimal models*, Mem. Soc. Math. France 9/10, 1983.

[20] ———, *Rational fibrations, minimal models, and fibrings of homogeneous spaces*, Trans. Amer. Math. Soc. **244** (1978) 199-224.

[21] P. Hilton, G. Mislin, and J. Roitberg, *Localization of nilpotent groups and spaces*, North-Holland Mathematics Studies, no. 15, New York, 1975.

[22] P. S. Hirschhorn, *Localization of model categories*, in preparation, available at http://www-math.mit.edu/~psh.

[23] M. J. Hopkins, *Formulations of cocategory and the iterated suspension*, Astérisque 113-114, Soc. Math. France (1984) 212-226.

[24] Mark Hovey, *Lusternik-Schnirelmann cocategory*, Illinois J. Math. **37** (1993) 224-239.

[25] ———, *Model categories*, Amer. Math. Soc., Mathematical surveys and monographs **63**, 1998.

[26] I. M. James, *Lusternik-Schnirelmann category*, in: Handbook of algebraic topolgoy, ed. I. M. James, Elsevier Science, 1995, 1293-1310.

[27] Daniel Quillen, *Homotopical algebra*, Lecture Notes in Math. **43**, Springer-Verlag, 1967.

[28] Mohammed Sbaï, *Cocatégorie rationnelle d'un espace topologique*, Astérisque 113-114 (1984) 288-291.

[29] ———, *La Cocatégorie Rationnelle d'un Espace*, Thèse de 3ème cycle, L'université de Lille, 1984.

[30] J.-P. Serre, *Groupes d'homotopie et classes de groupes abéliens*, Ann. Math. **58** (1953) 258-294.

[31] Dennis Sullivan, *Geometric topology, part I: Localization, periodicity, and Galois symmetry*, MIT mimeographed notes, 1970.

[32] K. Varadarajan, *Numerical invariants in homotopical algebra I, II*, Canad. J. Math. **27** (1975) 901-934, 935-960.

Department of Mathematics
Massachusetts Institute of Technology, 2-230
77 Massachusetts Avenue
Cambridge, MA 02139-4307
USA

E-mail: donald@math.mit.edu